\documentclass{article}
\usepackage{amssymb}
\usepackage{amsmath}
\usepackage{amsthm} 
\usepackage{mathtools}
\usepackage{bm} 
\usepackage{extarrows}
\usepackage{indentfirst}
\usepackage{stmaryrd}
\usepackage{tikz}
\usetikzlibrary{arrows, graphs}
\usepackage{setspace}
\def\ex{\mbox{ex}}

\def\im{\mbox{im}}

\newtheorem{thm}{Theorem}[section]
\newtheorem{cor}[thm]{Corollary}
\newtheorem{lem}[thm]{Lemma}

\newtheorem{claim}{Claim}[section]

\begin{document}
	  \renewcommand{\thefootnote}{\fnsymbol{footnote}}
	  \title{The cycle of length four is strictly $F$-Tur\'an-good
	  	\footnote{The work was supported by the National Natural Science Foundation of China (No. 12071453)
	  		and the National Key R and D Program of China(2020YFA0713100), the Anhui Initiative in
	  		Quantum Information Technologies (AHY150200) and the Innovation Program for Quantum
	  		Science and Technology, China (2021ZD0302904).}}
	  \author{Doudou Hei$^1$, Xinmin Hou$^{1,2}$\\
	  \small $^1$School of Mathematical Sciences\\
	  \small University of Science and Technology of China, Hefei, Anhui 230026, China\\
	  \small $^{2}$CAS Key Laboratory of Wu Wen-Tsun Mathematics\\
	  \small University of Science and Technology of China, Hefei, Anhui 230026, China
      }
      \date{}
      \maketitle

	   \begin{abstract}
	  Given an $(r+1)$-chromatic graph $F$ and a graph $H$ that does not contain $F$ as a subgraph, we say that $H$ is strictly $F$-Tur\'an-good if the Tur\'an graph $T_{r}(n)$ is the unique graph containing the maximum number of copies of $H$ among all $F$-free graphs on $n$ vertices for every $n$ large enough. 
	   Gy\H{o}ri, Pach and Simonovits (1991) proved that cycle $C_4$ of length four is strictly   $K_{r+1}$-Tur\'{a}n-good for all $r\geq 2$.  In this article, we extend this result and show that $C_4$ is strictly $F$-Tur\'an-good,  where  $F$ is an  $(r+1)$-chromatic graph with $r\ge 2$ and  a color-critical edge. Moreover, we show that every $n$-vertex $C_4$-free  graph $G$ with $N(H,G)=\ex(n,C_4,F)-o(n^4)$ can be obtained by adding or deleting $o(n^2)$ edges from $T_r(n)$. Our proof uses the flag algebra method developed by Razborov (2007).
   	   \end{abstract}

	  \section{Introduction}
	   All graphs considered in this article are finite and simple.  
	   Given a graph $G$, write $n(G)$ (resp. $e(G)$) for $|V(G)|$ (resp. $|E(G)|$) and write $G[Z]$ for the induced subgraph of $G$ on the vertex set $Z\subseteq V(G)$.  Let $X$ and $Y$ be disjoint subsets of $V(G)$. By $G[X,Y]$, we denote the bipartite subgraph of $G$ consisting of all edges that have one endpoint in $X$ and another in $Y$.
	   For mutually disjoint subsets $V_1, V_{2}, \ldots, V_{k}\subseteq V(G)$, similarly, we define $G[V_1, \ldots, V_{k}]$ to be the $k$-partite subgraph of $G$ consisting of all edges in ${\cup}_{1\leq i< j\leq k}E(G[V_i,V_j])$.
	   Write $K(V_1, \ldots, V_k)$ for the complete $k$-partite graph with color classes $V_1, \ldots, V_k$ and write $K_{t_1, \ldots, t_k}$ for a complete $k$-partite graph $K(V_1, \ldots, V_k)$ with $|V_i|=t_i$ for $i\in [k]$, where $[k]=\{1,2,\ldots,k\}$.


 Fix a graph $F$, we say that a graph $G$ is {\em $F$-free} (or {\em induced $F$-free}) if it does not contain $F$ as a subgraph (or an induced subgraph).	  
For given graphs $H$ and $G$, we define $N(H,G)$ (resp. $N_I(H,G)$) as the number of subgraphs (resp. induced subgraphs) of $G$ isomorphic to $H$. 
Let $\ex(n,H,F)$ denote the maximum value of $N(H,G)$ among all $F$-free graphs on $n$ vertices,
	  and we call the graph with $\ex(n, H, F)$ copies of $H$ an \emph{extremal graph}. This function is well-studied when $H$ is an edge, and it is called the Tur\'an number $\ex(n,F)$ of $F$ (one can eee, for example~\cite{Simonovits2013}, for a survey).
	  
 Let $T_r(n)$ denote the $r$-partite $n$-vertex Tur\'an graph, i.e., the  $r$-partite $n$-vertex complete graph of which each partite is of size $\lceil \frac{n}{r} \rceil$ or$ \lfloor \frac{n}{r} \rfloor$.  As pointed by Gerbner and Palmer~\cite{GERBNER2022103519}, there are few $F$-free graphs we have already known that are good candidates for being extremal constructions for maximizing copies of $H$, an exception is the Tur\'an graph, they call $H$ to be $F$-Tur\'an-good under this situation. More precisely, given an $(r+1)$-chromatic graph $F$ and a graph $H$ does not contain $F$ as a subgraph, we say that $H$ is {\it $F$-Tur\'{a}n-good} (or {\it strictly $F$-Tur\'an-good}) if $\ex(n,H,F)=N(H,T_{r}(n))$ (and the Tur\'an graph $T_{r}(n)$ is the unique extremal graph) for every $n$ large enough. We also call $(H,F)$ Tur\'an-good (or strictly Tur\'an-good) for short.

	  Gy\H{o}ri, Pach and Simonovits \cite{Ervin1991On} proved that $(C_4, K_{r+1})$ is strictly  Tur\'{a}n-good.	  
\begin{thm}[\cite{Ervin1991On}]\label{THM: GPS19}	
$C_4$ is strictly  $K_{r+1}$-Tur\'{a}n-good. Precisely, for every $K_{r+1}$-free graph $G$ with $|V(G)|=n\ge \max\{r, 5\}$, $N(C_4, G)\le N(C_4, T_{r}(n))$,  equality holds if and only if $G\cong T_{r}(n)$.
\end{thm}

 We say that an edge $e$ of a graph $F$ is \emph{color-critical} if deleting $e$ from $F$ results in a graph with a smaller chromatic number.
In extremal problems, a graph $F$ with $\chi(F)=r$ and a color-critical edge often behaves similarly to $K_r$. The most famous one along this flavor was given by Simonovits~\cite{S68}, who showed that $(K_2, F)$ is strictly $F$-Tur\'an-good, and Ma and Qiu~\cite{Ma} proved a generalized version by extending  the pair $(K_2, F)$ to $(K_r, F)$, where $\chi(F)>r\ge 2$.  There are also a few (strictly) Turán-good pairs given by researchers, including Gerbner, Palmer, Murphy, etc. Here is a list of some other Tur\'an-good pairs as we have known so far:

(1) (Gy\H{o}ri, Pach and Simonovits\cite{Ervin1991On}) $(H, K_{3})$ is strictly  Tur\'{a}n-good, where $H$ is a bipartite graph with matching number $\lfloor \frac{V(H)}{2}\rfloor$ (inclding the path $P_\ell$, the even cycle $C_{2\ell}$ and the Turán graph $T_2(m)$);

(2) (Gy\H{o}ri, Pach and Simonovits\cite{Ervin1991On}) $(K_{2,t}, K_{r})$ is strictly  Tur\'{a}n-good for $t=2,3$;

(3) (Gerbner and Palmer~\cite{GERBNER2022103519}) $(H, K_{k})$ is  Tur\'{a}n-good for $k\geq k_0$, where $H$ is a complete multipartite graph and $k_0$ is a consatant depending on $H$, and Gerbner and Palmer conjectured that this result is true for any graph $H$;

(4) (Gerbner~\cite{Gerbner4})  For any positive integers $m$ and $\ell$, $(P_m, C_{2\ell+1})$ and $(C_{2m}, C_{2\ell+1})$ are Tur\'{a}n-good. 

(5) (Gerbner and Palmer~\cite{GERBNER2022103519})   $(C_4, B_{2})$ and $(C_4, F_2)$ are Tur\'{a}n-good, where $B_k$ (resp. $F_k$) is the graph of $k$ triangles all sharing exactly one common edge (resp. one common vertex);

(6) (Gerbner and Palmer~\cite{Gerbner4})  For any positive integers $m$ and $t$, $(P_m, B_t)$ is Tur\'{a}n-good. 

(7)(Gerbner and Palmer~\cite{GERBNER2022103519}, Gerbner~\cite{Gerbner3})  $(P_3, F)$ is Tur\'{a}n-good, where $F$ is a  graph with $\chi(F)=k\geq 3$ and a color-critical edge;

(8)(Murphy and Nir~\cite{Murphy}, Qian et al~\cite{Qian})  $(P_4, K_{k})$ and $(P_5, K_{k})$ are Tur\'{a}n-good for $k\geq 4$;

(9)(Gerbner \cite{Gerbner3}) $(M_{2\ell}, F)$ is Tur\'{a}n-good, where $M_{2\ell}$ is a matching on $2\ell$ vertices,  $F$ is a  graph with $\chi(F)=k\geq 3$ and a color-critical edge.

Recently, a special family of bipartite graphs was proven to be strictly $F$-Tur\'an-good when $F$ is a graph with $\chi(F)=3$ and  a color-critical edge by Hei, Hou and Liu~\cite{unknown}. 
\begin{thm}[\cite{unknown}] \label{THM: F=3} 
	Let $F$ be a  graph with $\chi(F)=3$ and a color-critical edge and let $H$ be a bipartite graph with matching number $\left\lfloor \frac{|V(H)|}{2}\right\rfloor$.  Then $H$ is strictly $F$-Tur\'an good, i.e., $\ex(n, H, F)=N(H, T_2(n))$ for every $n$ large enough, and the Tur\'an graph $T_2(n)$ is the unique extremal graph for $(H, F)$.
\end{thm}
As a special case, we know that $(C_4, F)$ is strictly $F$-Tur\'an-good when $F$ is a graph with $\chi(F)=3$ and  a color-critical edge, which also has been shown in~\cite{Gebner-nonalign} (a special case of Theorem 1.9). 
\begin{cor}[\cite{Gebner-nonalign}]\label{COR: F=3}
	$(C_4, F)$ is  strictly $F$-Tur\'{a}n-good when $F$ is a graph with $\chi(F)=3$ and  a color-critical edge.
\end{cor}

In this article, we continue to show that Corollary~\ref{COR: F=3} also holds for $F$ with $\chi(F)=r+1\ge 3$ and a color-critical edge. This result is also a special case of results given by Gebner~\cite{Gebner22}, but that is essentially without proof, here we give the proof by a different way.  
\begin{thm}\label{THM: Main1}
	Let  $F$ be a graph with $\chi(F)=r+1\ge 4$ and  a color-critical edge. Then $(C_4, F)$ is strictly Tur\'{a}n-good.
\end{thm}

In fact,  Hei, Hou and Liu~\cite{unknown} have given a necessary and sufficient condition for a graph $H$ to be strictly $F$-Tur\'{a}n-good using the so called `weak $(r+1)$-T-property' and `T-extremal'.  We say a graph $H$ has the {\it weak $(r+1)$-T-property} if $N(H, K)\le N(H, T_r(n))$ for every complete $r$-partite graph $K=K_{t_1,\cdots,t_r}$ with $t_1+\cdots+t_r=n$ and the equality holds if and only if $K\cong T_r(n)$ for every $n$ large enough. Let $F$ be a graph with $\chi(F)=r+1$, an $n$-vertex $F$-free graph $G$ is called {\it T-extremal} if $\vert e(G)-e(T_r(n))\vert =o(n^2)$.
 \begin{thm}[\cite{unknown}]\label{THM: HHL} 
      	Let $F$ be a graph with $\chi(F)=r+1\ge3$ and a color-critical edge and let $H$ be a connected graph with $\chi(H)\le r$. Suppose every $n$-vertex $F$-free  graph $G$ with $N(H,G)=\ex(n,H,F)$ is T-extremal. If $H$ has the weak $(r+1)$-T-property, then $H$ is strictly $F$-Tur\'{a}n-good.
\end{thm}

\noindent{\bf Remark A:} Gerbner~\cite{Gebner22} defined a generalization of T-extremal $F$-free  graph $G$ with $N(H,G)=\ex(n,H,F)$: Let $\chi(H)<\chi(F)=k+1$. We say that $H$ is {\it $F$-Tur\'an-stable} if every $n$-vertex $F$-free graph $G$ with $N(H,G)\ge \ex(n,H,F)-o(n^{|V(H)|})$ can be
obtained from $T_k(n)$ by adding and removing $o(n^2)$ edges. In fact, alongside the proof of Theorem~\ref{THM: HHL} given in~\cite{unknown}, Theorem~\ref{THM: HHL} can be restated as follows: Let $F$ be a graph with $\chi(F)=r+1\ge3$ and a color-critical edge and let $H$ be a connected graph with $\chi(H)\le r$. Suppose $H$ is $F$-Tur\'an-stable. If $H$ has the weak $(r+1)$-T-property, then $H$ is strictly $F$-Tur\'{a}n-good.

\noindent{\bf Remark B:} By Theorem~\ref{THM: GPS19}, $C_4$ has the weak $(r+1)$-T-property. Therefore, by the above restated Theorem~\ref{THM: HHL}, to prove Theorem~\ref{THM: Main1},  it is sufficient to show that $C_4$ is $F$-Tur\'an-stable if $F$ is a graph with $\chi(F)=r+1\ge 4$ and  a color-critical edge.  
The following stability theorem for $C_4$ is another main result of this paper and has its own flavor in stability theory.
      
\begin{thm}\label{THM: Main2}
Let $F$ be a graph with $\chi(F)=r+1\ge 4$ and  a color-critical edge. Then $C_4$ is $F$-Tur\'an-stable.
\end{thm}

	
The rest of the article is arranged as follows.  In Section 2, we give the proof of Theorem~\ref{THM: Main2} admitting an important lemma (Lemma~\ref{LEM: inducedcherry2.6}). Section 3 will give a brief overview of the flag algebra. In the last section,  we prove Lemma~\ref{LEM: inducedcherry2.6} using the flag algebra method.

\section{ Proof of Theorem~\ref{THM: Main2}}

Let $H$ be a fixed graph. An $s$ $blow$-$up$ of a graph $H$ is the graph obtained by replacing each vertex $v$ of $H$ with an independent set $W_v$ of size $s$, and each edge $uv$ of $H$ with a complete bipartite graph between the corresponding two independent sets $W_u$ and $W_v$. We need the following nice results of graphs.
     
	 \begin{lem}[\cite{ALON2015683}]\label{LEM: blowup} Let $H$ be a fixed graph with $h$ vertices and let $F$ be a graph. Then $\ex ( n,H,F)$ $=\Omega (n^h)$ if and only if $F$ is not a subgraph of a blow-up of $H$. Otherwise, $\ex(n, H, F)\leq n^{h-\alpha}$ for some $\alpha>0$. 
	 \end{lem}

	 \begin{lem}[Induced Removal Lemma \cite{Alon}]\label{LEM: Remove} Let $\mathcal{F} $ be a set of graphs. For each $\varepsilon>0$, there exist $n_I(\varepsilon)>0$ and $\delta_I(\varepsilon)>0$ such that for every graph $G$ of order $n\geq n_I(\varepsilon)$, if $G$ contains at most $\delta_I(\varepsilon) n^{\vert V(H)\vert}$ induced copies of $H$ for every $H\in \mathcal{F} $, then $G$ can be made induced $\mathcal{F}$-free by removing or adding at most $\varepsilon n^2$ edges from $G$.
	 \end{lem}

	 \begin{lem}[\cite{Ma,Murphy}]\label{LEM: K_4} 
	 	Let $F$ be a graph with chromatic number $\chi(F)=r+1$ and  a color-critical edge. If $G$ is an $F$-free graph on $n$ vertices, then $N(K_4, G)\leq N(K_4, T_r(n))=\frac{r^3-6r^2+11r-6}{r^3}\binom{n}{4}+o(n^4).$
	 \end{lem}

	 \begin{lem}[\cite{Murphy}]\label{LEM: co-cherry} 
	Let  $P_3^c$ be the unique graph on three vertices with one edge (also called the co-cherry graph).
	 Then graph $G$ is a complete multipartite graph if and only if it does not contain the co-cherry graph $P_3^c$  as an induced subgraph.
	 \end{lem}

We first show the following lemma. 
	 \begin{lem}\label{LEM: 2.5} 
Let $f(x_1, x_2,\ldots, x_r)=N(C_4, K_{x_1,x_2,\ldots, x_r})$. If there is some $x_i\ge x_{j}+2$, then $f(x_1,\ldots, x_i-1,\ldots, x_j+1,\ldots,x_r)>f(x_1,\ldots, x_i,\ldots, x_j,\ldots,x_r)$.
	\end{lem}
 \begin{proof} We may assume $i=1, j=2$. Let $H=K(X_1, X_2, \ldots, X_r)$ with $|X_i|=x_i$ for $i\in[r]$ and let $H^*=K(X_1\setminus\{v\}, X_2\cup\{v^*\}, X_3,\ldots, X_r)$, where $v\in X_1$ and $v^*$ is a new  vertex added to $X_2$. 
	 Denote $X_1^*=X_1\setminus\{v\}$ and $X_2^*=X_2\cup\{v^*\}$. 
	 Let $N_H(v, C_4)=\{C : C\cong C_4 \text{ in $H$ with } v\in V(C)\}$ and $N_{H^*}(v^*, C_4)=\{C : C\cong C_4 \text{ in $H^*$ with } v^*\in V(C)\}$.
	 Let $n_H(v, C_4)=|N_H(v, C_4)|$ and $n_{H^*}(v^*, C_4)=|N_{H^*}(v^*, C_4)|$. Since the copy of $C_4$ that does not pass through $v$ in $H$ remains unchanged in $H^*$, to show the lemma, it suffices to show that $n_H(v, C_4)<n_{H^*}(v^*, C_4)$. 
	In addition, a copy $C\in N_H(v, C_4)$ with $V(C)\cap(X_1\cup X_2)=\{v\}$ corresponds to a copy $C^*\in N_{H^*}(v^*, C_4)$ with $V(C^*)\cap (X_1^*\cup X_2^*)=\{v^*\}$ and vice versa. 
It suffices to focus on those $C_4$ that contain $v$ (or $v^*$) and at least one other vertex in $X_1\cup X_2$ (or in $X_1^*\cup X_2^*$). 
   Let $$c(v, n_1, n_2)=|\{C\in N_H(v, C_4) : |V(C)\cap X_i|=n_i \text{ for } i=1,2\}|$$ and $$c^*(v^*, n_1^*, n_2^*)=|\{C\in N_{H^*}(v^*, C_4) : |V(C)\cap X_i^*|=n_i^* \text{ for } i=1,2\}|.$$
   Then we have  $1\le n_1\le 2$, $0\le n_2\le 2$, $0\le n_1^*\le 2$, $1\le n_2^*\le 2$, and $c(v, 1,0)=c^*(v^*, 0,1)$. 
Therefore, we only need to show $c(v,1,1)+c(v,1,2)+c(v,2,0)+c(v,2,1)+c(v,2,2)<c^*(v^*,1,1)+c^*(v^*,2,1)+c^*(v^*,0,2)+c^*(v^*,1,2)+c^*(v^*,2,2)$.
Let $I=\{3,4,\cdots,r\}.$
	 We count the number of $C_4$ according to the choices of $n_1$ and $n_2$,
\begin{eqnarray*}
&&	c(v,1,1)=x_2\cdot\left(\sum_{i\in I}\binom{x_i}{2}+3\sum_{\{i,j\}\in \binom{I}{2}}x_ix_j\right),\\
&&	c(v,1,2)=\binom{x_2}{2}\sum_{i\in I}x_i,          \\
&& c(v,2,0)=(x_1-1)\left(\sum_{i\in I}\binom{x_i}{2}+\sum_{\{i,j\}\in \binom{I}{2}}x_ix_j\right),\\
&& c(v,2,1)=(x_1-1)x_2(n-x_1-x_2),\\
&& c(v,2,2)=(x_1-1)\binom{x_2}{2}.
\end{eqnarray*}	
Similarly,	we have 
\begin{eqnarray*}
&&c^*(v^*,1,1)=(x_1-1)\cdot\left(\sum_{i\in I}\binom{x_i}{2}+3\sum_{\{i,j\}\in \binom{I}{2}}x_ix_j\right),\\
&&c^*(v^*,2,1)=\binom{x_1-1}{2}\sum_{i\in I}x_i,\\
&&c^*(v^*,0,2)=x_2\cdot\left(\sum_{i\in I}\binom{x_i}{2}+\sum_{\{i,j\}\in \binom{I}{2}}x_ix_j\right),\\
&& c^*(v^*,1,2)=(x_1-1)x_2(n-x_1-x_2),\\
&&c^*(v^*,2,2)=x_2\binom{x_1-1}{2}.
\end{eqnarray*}
Since $x_1\geq x_2+2$, we have $c(v,1,2)<c^*(v^*,2,1)$, $c(v,2,1)=c^*(v^*,1,2)$, $c(v,2,2)<c^*(v^*,2,2)$, and 
\begin{eqnarray*}
	c(v,1,1)+c(v,2,0)&=&(x_1+x_2-1)\cdot\sum_{i\in I}\binom{x_i}{2}+(x_1+3x_2-1)\sum_{\{i,j\}\in \binom{I}{2}}x_ix_j\\
&<& (x_1+x_2-1)\cdot\sum_{i\in I}\binom{x_i}{2}+(3x_1+x_2-3)\sum_{\{i,j\}\in \binom{I}{2}}x_ix_j\\
&=& c^*(v^*,1,1)+c^*(v^*,0,2).
\end{eqnarray*}
This completes the proof. 
\end{proof}

	   
   

     Let $H$ and $G$ be graphs on $n_1$ and $n_2$ vertices, respectively, where $n_1\leq n_2$. The $density$ $d(H,G)$ of $H$ in $G$ is defined by 
     $$d(H,G)=\frac{N(H,G)}{\binom{n_2}{n_1}}.$$ 
     Let $F$ be a graph with $\chi(F)=r+1$ and a color-critical edge, and let $\mathcal{F}_{n,r}$ denote the family of $F$-free graphs on $n$ vertices.
      Define $OPT_r(C_4)$ as follows:
     $$OPT_r(C_4)=\lim_{n \rightarrow \infty} \max_{G\in \mathcal{F}_{n,r} }d(C_4, G).$$
     
\begin{lem}\label{LEM: inducedcherry2.6}  
	Let $F$ be a graph with $\chi(F)=r+1\ge 4$ and a color-critical edge.
	For any $\delta>0$, there exist $n_{c}=n_c(\delta) $ and $\varepsilon_{c}=\varepsilon_c(\delta)>0$ such that for every $F$-free graph $G$ of order $n\geq n_{c}$, if $d(C_4, G)\geq OPT_r(C_4)-\varepsilon_{c}$,  then $G$ contains at most $\delta n^3$ induced copies of the co-cherry graph $P_3^c$.
\end{lem}

The above lemma plays an important role in the proof of Theorem~\ref{THM: Main2}, and we will prove it in Section 4 by the flag algebra argument.
The following is a stability lemma for almost optimal complete $r$-partite graphs. 
\begin{lem}\label{LEM: stability}
Let $G$ be a complete $r$-partite graph with partite sets $X_1, \cdots, X_r$. For any $\varepsilon>0$, there exists $\delta_s=\delta_s(\varepsilon)>0$ such that if 
$d(C_4, G)>OPT_r(C_4)-\delta_s,$
then for each $i=1,2,\cdots,r$,
$$\frac{(1-\varepsilon)|V(G)|}{r}\leq \vert X_i \vert \leq \frac{(1+\varepsilon)|V(G)|}{r}.$$
\end{lem}
\begin{proof}
\indent Let $|V(G)|=n$ and let $x_i=\vert X_i \vert$ for $i=1,2,\ldots, r$.  Without loss of generality, we may assume that $x_1=\frac{1+\eta(r-1)}{r}n>\frac{n(1+\varepsilon)}{r}$. Then  $\eta(r-1)>\varepsilon$. 
By Lemma~\ref{LEM: 2.5}, $d(C_4, {G})$ is maximized if all the remaining parts are balanced, i.e. $x_i=\frac{(1-\eta)n}{r}$ for  $i= 2,\ldots,r$. Let ${G}^{*}$ be the graph induced by the vertex set $V({G})\setminus X_1$. Then
$$N(C_4, {G}^*)=N\left(C_4, T_{r-1}\left(\frac{(r-1)(1-\eta)n}{r}\right)\right)+o(n^4).$$
So
\begin{equation*}
	\begin{aligned}
		N(C_4,G)
		&= N(C_4, G^*)+3x_1\sum_{2\le i<j<k\le r}x_ix_jx_k+2x_1\sum_{2\le i<j\le r}\binom{x_i}{2}x_j\\
		&\quad +\binom{x_1}{2}\sum_{2\le i<j\le r}x_ix_j+\binom{x_1}{2}\sum_{i=2}^r\binom{x_i}{2}\\
		&=N\left(C_4, T_{r-1}\left(\frac{(r-1)(1-\eta)n}{r}\right)\right)\\
		&\quad +\frac{(r-1)(r-2)^2(1+\eta(r-1))(1-\eta)^3}{2r^4} n^4\\
		&\quad +\frac{(r-1)^2 (1+\eta(r-1))^2(1-\eta)^2}{4r^4} n^4+o(n^4).
	\end{aligned}
\end{equation*}
Therefore,
\begin{equation}
	\begin{aligned}
		d(C_4,{G})&=\frac{N(C_4,{G})}{\binom{n}{4}}\\
		&=OPT_{r-1}(C_4)\left(\frac{1-\eta}{r}\right)^4+\frac{12(r-1)(r-1)^2(1+\eta(r-1))(1-\eta)^3}{r^4}\\
		&\quad +\frac{6(r-1)^2(1+\eta(r-1))^2(1-\eta)^2}{r^4}+o(1)\\
		&{=OPT_r(C_4)-6\eta^2(\frac{2r^3-10r^2+17r-9}{r^3})+12\eta^3\frac{r^3 - 6r^2 + 11r - 6}{r^3}}\\
		&{\quad -3\eta^4 \frac{r^3 - 8r^2 + 16r - 9}{r^3}+o(1).}
	\end{aligned}
	\nonumber
\end{equation}
Let $g(r)=\frac{2r^3-10r^2+17r-9}{r^3}=2-\frac{10}{r}+\frac{17}{r^2}-\frac{9}{r^3}$, which is positive with minimum $g(4)=\frac{27}{64}$. Therefore, for sufficently small $\eta$ and large $n$, we get $ d(C_4,\mathcal{G})\leq OPT_r(C_4)-2\eta^2, $ a contradiction. This implies the statement of the claim. 
\end{proof}
     
 As pointed out in Remark A, Theorem~\ref{THM: Main2} implies Theorem~\ref{THM: Main1}. Therefore, it is sufficient to show Theorem~\ref{THM: Main2}, i.e., we will show that every $n$-vertex $C_4$-free  graph $G$ with $N(H,G)=\ex(n,C_4,F)-o(n^4)$ can be obtained by adding or deleting $o(n^2)$ edges from $T_r(n)$.  
\begin{proof}[Proof of Theorem~\refeq{THM: Main2}]
Let $F$ be a graph with chromatic number $\chi(F)=r+1$ and  a color-critical edge. Let $G$ be an $n$-vertex $F$-free graph with $N(H,G)=\ex(n,C_4,F)-o(n^4)$.

First, let us consider the case $r= 3$. 
Since $F$ is a subgraph of a blow-up of $K_4$, by Lemma~\ref{LEM: blowup}, 
$$N(K_4, G)\leq \ex(n, K_4, F)=o(n^4).$$
By Lemma~\ref{LEM: Remove}, we can remove or add  $o(n^2)$ edges from $G$ and obtain a $K_4$-free graph $G^*$. The removal of $o(n^2)$ edges from $G$ can destroy at most $o(n^2)\cdot O(n^{2})=o(n^4)$ copies of $C_4$. Therefore, it suffices to estimate the number of $C_4$  and the number of edges in $G^*$.
Let $$\mathcal{M}_2=\{M_2 : M_2=\{e_1, e_2\} \text{ is an independent set in } G^*\}.$$
Since $G^*$ is $K_4$-free, the number of $C_4$ contained in the induced graph of a fixed $M_2\in\mathcal{M}_2$ is at most 1. Note that each copy of $C_4$ in $G^*$ contains exactly 2 members in $\mathcal{M}_2$. Therefore, we have  
$$N(C_4, T_3(n)) -o(n^4)\le N(C_4, G)-o(n^4)\le N(C_4, G^*)\le\frac 12|\mathcal{M}_2|.$$
Since $N(C_4, T_3(n))=\frac {n^4}{36}+o(n^4)$ by (\ref{EQ: NumC_4inTuran}) and  $\frac 12|\mathcal{M}_2|\le\frac{1}{2}\binom{e(G^*)}{2}$, we have $$e(G)\geq e(G^*)-o(n^2)\geq \frac{1}{3}n^2-o(n^2)=e(T_3(n))-o(n^2).$$

So from now on we assume that $r\geq 4$. 
It suffices to prove that for any $\varepsilon>0$, there exist $n_0>0$ and $\varepsilon'>0$ such that for every $F$-free graph $G$ of order $n\geq n_0$, if $d(C_4, G)\geq OPT_r(C_4)-\varepsilon'$,  then by changing at most $\varepsilon n^2$ pairs of adjacencies, we can obtain $T_r(n)$ from $G$.

{Let $\mathcal{F}=\{K_{r+1}, P_3^c\}$ and $\varepsilon_1>0$. By Lemma~\ref{LEM: Remove}, we have $\delta_I(\varepsilon_1)>0$ and $n_I(\varepsilon_1)>0$. For $\delta_I(\varepsilon_1)>0$,
	since $F$ is a subgraph of a blow-up of $K_{r+1}$, by Lemma~\ref{LEM: blowup}, there exists an integer $n_F(\delta_I(\varepsilon_1))>0$, such that every $F$-free graph on $n>n_F(\delta_I(\varepsilon_1))$ vertices contains at most $\delta_I(\varepsilon_1) n^{r+1}$ copies of $K_{r+1}$. 
	By Lemma~\ref{LEM: inducedcherry2.6}, we also have integer $n_{c}(\delta_I(\varepsilon_1))>0$ and $\varepsilon_{c}(\delta_I(\varepsilon_1))>0$ such that for every $F$-free graph $G$ of order $n\geq n_c(\delta_I(\varepsilon_1))$, if $d(C_4, G)\geq OPT_r(C_4)-\varepsilon_{c}(\delta_I(\varepsilon_1))$,  then $G$ contains at most $\delta_I(\varepsilon_1) n^3$ induced copies of the co-cherry graph $P_3^c$.
Now we choose $n_0\geq \max\{n_F(\delta_I(\varepsilon_1)), n_{c}(\delta_I(\varepsilon_1)), n_I(\varepsilon_1)\}$ and $\varepsilon_2<\varepsilon_{c}(\delta_I(\varepsilon_1))$. 
	Assume that $G$ is a graph on $n\geq n_0$ vertices such that 
	$$d(C_4,G)>OPT_r(C_4)-\varepsilon_2.$$
Therefore, Lemmas~\ref{LEM: blowup} and~\ref{LEM: inducedcherry2.6} guarantee that $G$ contains at most $\delta_I(\varepsilon_1) n^{r+1}$ copies of $K_{r+1}$ and at most $\delta_I(\varepsilon_1) n^{3}$ induced copies of $P_3^c$.}
By Lemma~\ref{LEM: Remove}, we can make $G$ into an  induced $\mathcal{F}$-free graph $G^*$ by deleting or adding at most $\varepsilon_1 n^2$ edges.  Since one removed edge from $G$ destroys at most $n^2$ copies of $C_4$, the total number of destroyed copies of $C_4$ is at most $\varepsilon_1n^2\cdot n^2=\varepsilon_1n^4$.
So $N(C_4, G^*)\geq N(C_4, G)-\varepsilon_1n^4$. 
Therefore, for $n$ large enough, we have
$$d(C_4,G^*)\geq OPT_r(C_4)-24\varepsilon_1-\varepsilon_2,$$
Since $G^*$ is induced $P_3^c$-free, by Lemma~\ref{LEM: co-cherry}, $G^*$ is complete multipartite. Since $G^*$ is $K_{r+1}$-free, we may assume $G^*$ is a complete $r$-partite graph with partite sets $X_1, \cdots, X_r$. We will complete the proof by showing that the partite sets are almost blanced.
Now,  we apply Lemma~\ref{LEM: stability} to $G^*$, for any $\varepsilon>0$, choose $\varepsilon_1\leq \frac{\varepsilon}{2}$ small enough and  $\delta_s(\frac{\varepsilon}{2})>24\varepsilon_1+\varepsilon_2$.
Then $d(C_4,G^*)\geq OPT_r(C_4)-24\varepsilon_1-\varepsilon_2 > OPT_r(C_4)-\delta_s(\frac{\varepsilon}{2})$. By Lemma~\ref{LEM: stability}, 
we have $$\frac{n(1-{\varepsilon}/{2})}{r}\leq \vert X_i \vert \leq \frac{n(1+{\varepsilon}/{2})}{r}$$ for $1\le i\le r$.
Therefore, by changing at most $(\varepsilon_1+ {\varepsilon}/{2})n^2 \leq \varepsilon n^2$ edges, we can obtain $T_r(n)$ from the original graph $G$, which completes the proof. 
\end{proof}	 
	   
\section{Overview of flag algebra}
 In this section, we give a brief overview of the  flag algebra method developed by Razborov \cite{Razborov}, which provides a frame work for computationally solving problems in extremal combinatorics. By this method, we can find inequalities of subgraph densities in graph limits with the help of semi-definite programming. First, let us present a brief introduction and description of the notation and theory needed in this section. 
 
 Let $H$ and $G$ be graphs on $n_1$ and $n_2$ vertices , respectively, where $n_1\leq n_2$. Recall that the $density$ $d(H,G)$ of $H$ in $G$ is defined by 
 $$d(H,G)=\frac{N(H,G)}{\binom{n_2}{n_1}}.$$ 
 Similarly, the {\it induced density} $P(H,G)$ of $H$ in $G$ is defined by
 $$P(H,G)=\frac{N_I(H,G)}{\binom{n_2}{n_1}}.$$ 
 Let a subset $Y$ be selected uniformly at random from $V(G)$ such that $\vert Y\vert=n_1$. Then we can interpret the {\it induced density} $P(H,G)$ of $H$ in $G$ as the probability that $G[Y]$ is isomorphic to $H$.
 
We will need the notion of a flag. A {\em type} of size $k$ is a graph $\sigma $ on $k$ vertices labeled by $[k]$. If $\sigma $ is a type of size $k$ and $F$ is a graph on at least $k$ vertices, then an {\em embedding} of $\sigma$ into $F$ is an injective function $\theta :[k]\mapsto V(F)$ such that $\theta $ gives an isomorphism between $\sigma $ and $F[\im(\theta )]$. A {\em $\sigma $-flag} is a pair $(F,\theta)$ where $F$ is a graph and $\theta$ is an embedding of $\sigma$ into $F$. Two $\sigma$-flags $(F,{\theta}_1)$ and $(G,{\theta}_2)$ are {\em isomorphic} if there exists a graph isomorphism between $F$ and $G$ that preserves the labeled subgraph $\sigma $.

Let $\mathcal{F} $ (or ${\mathcal{F} }_\ell$) denote the set of all graphs (or all graphs on $\ell$ vertices) up to isomorphism. Let ${\mathcal{F}}^\sigma$ (or ${\mathcal{F}}^\sigma_\ell $) denote the set of all $\sigma $-flags (or on $\ell$ vertices). Let $\mathbb{R}\mathcal{F}$, $\mathbb{R}\mathcal{F}^\sigma$ and $\mathbb{R}\mathcal{F}^\sigma_\ell$ denote the set of all formal finite linear combinations of elements in $\mathcal{F}$, ${\mathcal{F}}^\sigma$ and ${\mathcal{F}}^\sigma_\ell$, respectively, where the coefficients are real numbers. Note that if the size of $\sigma $ is 0, then ${\mathcal{F}}^\sigma=\mathcal{F}$  and ${\mathcal{F}}^\sigma_\ell=\mathcal{F}_\ell$.
 
 For two $\sigma$-flags $(H,\vartheta )$ and $(G,{\theta})$ with $n(H)\leq n(G)$, let $P((H,{\vartheta }),(G,{\theta}))$ denote the probability that any injective map from $V(H)$ to $V(G)$ that fixes the labeled graph $\sigma$ induces a copy of $H$ in $G$. Observe that if $\sigma $ is the empty graph, then $P((H,{\vartheta }),(G,{\theta}))=P(H,G)$.
 
 Let $(F_1, \theta _1)$,  $(F_2, \theta _2)$, $(G,\theta)\in {\mathcal{F}}^\sigma$ be three $\sigma$-flags for which $n(F_1)+n(F_2)\leq n(G)+n(\sigma)$.
 Let $X_1$ and $X_2$ be two disjoint sets of sizes $n(F_1)-n(\sigma)$ and $n(F_2)-n(\sigma)$ respectively, selected uniformly at random from $V(G)\setminus \im(\theta)$. Let $P((F_1, \theta _1),(F_2, \theta _2);(G, \theta ))$ denote the probability that $(G[X_1\cup \im(\theta )], \theta)$ is isomorphic to $(F_1, \theta_1 )$ and $(G[X_2\cup \im(\theta )], \theta)$ is isomorphic to $(F_2, \theta_2 )$. Razborov showed that as $n(G)$ grows,  the following inequality holds:
 $$\vert P((F_1, \theta _1),(F_2, \theta _2);(G, \theta ))-P((F_1,{\theta}_1),(G,{\theta}))P((F_2,{\theta}_2),(G,{\theta}))\vert \leq O({n(G)}^{-1}).$$
 Hence, as the size of $G$ tends to infinity, we can assume that we select $X_1$ and $X_2$ independently.
 
 Let $\mathcal{K}^\sigma $ denote the linear subspace of $\mathbb{R} \mathcal{F}^\sigma $ generated by all elements of the form
 $$(H,{\vartheta })-\sum_{(G,{\theta})\in {\mathcal{F}}^\sigma _m}P((H,{\vartheta }),(G,{\theta}))\cdot (G, \theta)$$
 where $m>n(H)$.
 Razborov has shown that there exists an algebra ${\mathcal{A}^\sigma }=\mathbb{R} \mathcal{F}^\sigma/\mathcal{K}^\sigma $ with well defined  addition and multiplication. Addition is defined in the natural way, by simply adding the coefficients of the elements in $\mathbb{R} \mathcal{F}^\sigma $. 
 Let $w=n(F_1)+n(F_2)-n(\sigma )$, then the product of $(F_1, \theta_1)$ and $(F_2, \theta_2)$ is defined as 
 $$(F_1, \theta_1)\cdot (F_2, \theta_2)=\sum_{(F_3, \theta_3)\in \mathcal{F}^\sigma_w} P((F_1, \theta _1),(F_2, \theta _2);(F_3, \theta_3 ))\cdot (F_3, \theta_3 ).$$ Addition and multiplication in ${\mathcal{A}^\sigma }$ are defined as an extension of addition and multiplication in $\mathbb{R}\mathcal{F}^\sigma$, respectively. If the size of $\sigma$ is 0, then we use $\mathcal{A}$ to denote  ${\mathcal{A}^\sigma }$.
 
 A sequence of graphs $\left(G_{n}\right)_{n \geq 1}$, where $n(G_{n})=n$, is said to be {\it convergent} if for every finite graph $H$, the $\lim\limits_{n \rightarrow \infty} P\left(H, G_{n}\right)$ exists. Let $\lim\limits_{n \rightarrow \infty}P( * , G_n)$ denote the corresponding linear function from $\mathbb{R}\mathcal{F}$ to $\mathbb{R}$.
 Let $\operatorname{Hom}^{+}\left(\mathcal{A}, \mathbb{R}\right)$ denote the set of all homomorphisms $\phi$ from $\mathcal{A}$ to $\mathbb{R}$ such that $\phi(F) \geq 0$ for each element $F \in \mathcal{F}$. Razborov showed that each function $\phi \in \operatorname{Hom}^{+}\left(\mathcal{A}, \mathbb{R}\right)$ corresponds to some convergent graph sequence $\left(G_{n}\right)_{n \geq 1}$, specifically, we have the following theorem. 
{\begin{thm}[\cite{Razborov}]
     (a) For every convergent sequence  $\left(G_{n}\right)_{n \geq 1}$,
     $$\lim_{n \rightarrow \infty}P( * , G_n)\in \operatorname{Hom}^{+}\left(\mathcal{A}, \mathbb{R}\right). $$
     (b) Conversely, every element $\operatorname{Hom}^{+}\left(\mathcal{A}, \mathbb{R}\right)$ can be represented in the form $$\lim_{n \rightarrow \infty}P( *, G_n)$$ for a convergent sequence $\left(G_{n}\right)_{n \geq 1}$.
 \end{thm}}
 For each type $\sigma$ labeled by $[k]$, Razborov also defined an unlabeling operator $$\llbracket \ \rrbracket_{\sigma}:  \mathcal{A}^{\sigma} \rightarrow  \mathcal{A}.$$  For a $\sigma$-flag $(F, \theta)$, let $q_{\sigma}(F)$ denote the probability that $\left(F, \theta^{\prime}\right)$ is isomorphic to $(F, \theta)$, where $\theta^{\prime}: [k] \rightarrow V(F)$ is a randomly chosen injective mapping. Let $F^{\prime}$ denote the graph isomorphic to $F$ when ignoring labels. Then
 $$\llbracket (F, \theta) \rrbracket_{\sigma}=q_{\sigma}(F) F^{\prime} .$$
 In addition, Razborov proved the following useful inequality.
 
 \begin{thm}[\cite{Razborov}]
 	Let $\sigma$ be a type and $\phi \in \operatorname{Hom}^{+}(\mathcal{A}, \mathbb{R})$, then, for any  $\alpha \in \mathcal{A}^{\sigma}$,
 	\begin{equation}\label{EQ: sigma}
 		\phi\left(\llbracket \alpha\cdot \alpha \rrbracket_{\sigma}\right) \geq 0.
 	\end{equation}
 \end{thm}
 \noindent Let $\left(G_{n}\right)_{n \geq 1}$ be the corresponding convergent sequence of $\phi$, then, by the above theorem, we have for any  $\alpha \in  \mathcal{A}^{\sigma}$,
 \begin{equation}\label{EQ: sigma2}
 	\lim_{n \rightarrow \infty}P\left(\llbracket \alpha\cdot\alpha \rrbracket_{\sigma}, G_n\right) \geq 0.
 \end{equation}

\section{Proof of Lemma~\ref{LEM: inducedcherry2.6}} 
 Now we are ready to give the proof of Lemma \ref{LEM: inducedcherry2.6}. 
 \begin{lem}[Restatement of Lemma 2.6]\label{Restate of Lemma 2.6}  
 	Let $F$ be a graph with $\chi(F)=r+1\ge 4$ and a color-critical edge.
 	For any $\delta>0$, there exist $n_{c}=n_c(\delta) $ and $\varepsilon_{c}=\varepsilon_c(\delta)>0$ such that for every $F$-free graph $G$ of order $n\geq n_{c}$, if $d(C_4, G)\geq OPT_r(C_4)-\varepsilon_{c}$,  then $G$ contains at most $\delta n^3$ induced copies of the co-cherry graph $P_3^c$.
 \end{lem}
 The outline of the proof is as follows: the main step is to calculate $OPT_r(C_4)$. We first use the number of $C_4$ in $T_r(n)$ to give a lower bound and then use the flag algebra argument to show that this lower bound is also an upper bound. To show the upper bound, we will calculate four functions $Q_i(r)$ and $q_i(r)$ for $i\in\{0,1,2,3\}$, which is inspired by \cite{Bernard}. Using semidefinite programming (interested readers can refer to \cite{Grzesik}), we can verify that the upper bound of $OPT_r(C_4)$ is correct for small values of $r$. After doing so, we are able to guess the prospective types and the order of the corresponding flags. Since the number of flags on 3 vertices is very small, we can greedily determine which specific flags to use with the help of MATLAB. 
 \begin{proof}
 	 First,  we give the following claim.    
 	\begin{claim}\label{CLM: claim 1}
 		For all $r\geq 3$, $OPT_r(C_4)=\frac{3(r-1)(r^2-3r+3)}{r^3}.$
 	\end{claim}
 	\begin{proof}[Proof of Claim \ref{CLM: claim 1}:] 
 		By definition, $$OPT_r(C_4)\ge \lim_{n \rightarrow \infty}\frac{N(C_4, T_r(n))}{{n\choose 4}}.$$ 
 		First we count $N(C_4, T_r(n))$ according to the distribution of $V(C_4)$ in $T_r(n)$. Let $N_i$ be the number of copies of $C_4$ with $V(C_4)$ distributed in exactly $i$ classes of $T_r(n)$, where $i=2, 3$ or $4$. So 
 		$$N_2=\binom{r}{2}\cdot \binom{\frac{n}{r}}{2}\cdot \binom{\frac{n}{r}}{2}+o(n^4)=\frac{(r-1)n^2(n-r)^2}{8r^3}+o(n^4),$$
 		$$N_3=\binom{r}{3}\cdot 3 \binom{\frac{n}{r}}{2}\cdot \frac{n}{r}\cdot \frac{n}{r}+o(n^4)=\frac{(r-1)(r-2)n^3(n-r)}{4r^3}+o(n^4),$$
 		and 
 		$$N_4=\binom{r}{4}\cdot {\left(\frac{n}{r}\right)}^4\cdot 3+o(n^4)=\frac{(r-1)(r-2)(r-3)n^4}{8r^3}+o(n^4),$$
 		where the error term $o(n^4)$ accounts for the cases when $n$ is not divisible by $r$.
 		Thus the total number of copies of $C_4$ in $T_r(n)$ is 
 		\begin{equation}\label{EQ: NumC_4inTuran}
 			N(C_4, T_r(n))=N_2+N_3+N_4=\frac{(r-1)(r^2-3r+3)n^4}{8r^3}+o(n^4).
 		\end{equation}     
 		Therefore,
 		\begin{equation}\label{EQ: lower}
 			OPT_r(C_4)\geq \lim_{n \rightarrow \infty} \frac{N(C_4, T_r(n))}{\binom{n}{4}}=\frac{3(r-1)(r^2-3r+3)}{r^3}.
 		\end{equation}
 		
 		To complete the proof of the claim, we show that the lower bound given in (\ref{EQ: lower})  is also an upper bound of $OPT_r(C_4)$. Let $\mathcal{F}_4=\{F_0, \ldots, F_{10}\}$ be the set of unlabeled graphs up to isomorphism on four vertices as drawn in Fig.~\ref{PIC: f4}.
 		
 		\begin{spacing}{2.0}\label{PIC: f4}
 			\qquad\qquad $F_0=$
 			\begin{tikzpicture}[baseline=10pt,scale=0.75]
 				\filldraw[fill=black,draw=black] (0,0) circle (0.05);
 				\filldraw[fill=black,draw=black] (0,1) circle (0.05);
 				\filldraw[fill=black,draw=black] (1,0) circle (0.05);
 				\filldraw[fill=black,draw=black] (1,1) circle (0.05);
 			\end{tikzpicture}
 			\quad $F_1=$
 			\begin{tikzpicture}[baseline=10pt,scale=0.75]
 				\draw[style=thick,color=blue] (0,0) -- (1,0);
 				\filldraw[fill=black,draw=black] (0,0) circle (0.05);
 				\filldraw[fill=black,draw=black] (0,1) circle (0.05);
 				\filldraw[fill=black,draw=black] (1,0) circle (0.05);
 				\filldraw[fill=black,draw=black] (1,1) circle (0.05);
 			\end{tikzpicture}
 			\quad $F_2=$
 			\begin{tikzpicture}[baseline=10pt,scale=0.75]
 				\draw[style=thick,color=blue] (0,0) -- (1,0);
 				\draw[style=thick,color=blue] (0,1) -- (0,0);
 				\filldraw[fill=black,draw=black] (0,0) circle (0.05);
 				\filldraw[fill=black,draw=black] (0,1) circle (0.05);
 				\filldraw[fill=black,draw=black] (1,0) circle (0.05);
 				\filldraw[fill=black,draw=black] (1,1) circle (0.05);
 			\end{tikzpicture}
 			\quad $F_3=$
 			\begin{tikzpicture}[baseline=10pt,scale=0.75]
 				\draw[style=thick,color=blue] (0,0) -- (1,0);
 				\draw[style=thick,color=blue] (0,1) -- (0,0);
 				\draw[style=thick,color=blue] (0,0) -- (1,1);
 				\filldraw[fill=black,draw=black] (0,0) circle (0.05);
 				\filldraw[fill=black,draw=black] (0,1) circle (0.05);
 				\filldraw[fill=black,draw=black] (1,0) circle (0.05);
 				\filldraw[fill=black,draw=black] (1,1) circle (0.05);
 			\end{tikzpicture}\\
 			
 			\qquad\qquad$F_4=$
 			\begin{tikzpicture}[baseline=10pt,scale=0.75]
 				\draw[style=thick,color=blue] (1,0) -- (1,1);
 				\draw[style=thick,color=blue] (1,1) -- (0,1);
 				\draw[style=thick,color=blue] (0,1) -- (1,0);
 				\filldraw[fill=black,draw=black] (0,0) circle (0.05);
 				\filldraw[fill=black,draw=black] (0,1) circle (0.05);
 				\filldraw[fill=black,draw=black] (1,0) circle (0.05);
 				\filldraw[fill=black,draw=black] (1,1) circle (0.05);
 			\end{tikzpicture}
 			\quad $F_5=$
 			\begin{tikzpicture}[baseline=10pt,scale=0.75]
 				\draw[style=thick,color=blue] (0,0) -- (1,0);
 				\draw[style=thick,color=blue] (1,1) -- (0,1);
 				\filldraw[fill=black,draw=black] (0,0) circle (0.05);
 				\filldraw[fill=black,draw=black] (0,1) circle (0.05);
 				\filldraw[fill=black,draw=black] (1,0) circle (0.05);
 				\filldraw[fill=black,draw=black] (1,1) circle (0.05);
 			\end{tikzpicture}
 			\quad $F_6=$
 			\begin{tikzpicture}[baseline=10pt,scale=0.75]
 				\draw[style=thick,color=blue] (0,0) -- (1,0);
 				\draw[style=thick,color=blue] (1,0) -- (1,1);
 				\draw[style=thick,color=blue] (0,1) -- (0,0);
 				\draw[style=thick,color=blue] (0,0) -- (1,1);
 				\filldraw[fill=black,draw=black] (0,0) circle (0.05);
 				\filldraw[fill=black,draw=black] (0,1) circle (0.05);
 				\filldraw[fill=black,draw=black] (1,0) circle (0.05);
 				\filldraw[fill=black,draw=black] (1,1) circle (0.05);
 			\end{tikzpicture}
 			\quad $F_7=$
 			\begin{tikzpicture}[baseline=10pt,scale=0.75]
 				\draw[style=thick,color=blue] (0,0) -- (1,0);
 				\draw[style=thick,color=blue] (1,1) -- (0,1);
 				\draw[style=thick,color=blue] (0,1) -- (0,0);
 				\filldraw[fill=black,draw=black] (0,0) circle (0.05);
 				\filldraw[fill=black,draw=black] (0,1) circle (0.05);
 				\filldraw[fill=black,draw=black] (1,0) circle (0.05);
 				\filldraw[fill=black,draw=black] (1,1) circle (0.05);
 			\end{tikzpicture}\\
 			
 			\qquad\qquad\qquad$F_8=$
 			\begin{tikzpicture}[baseline=10pt,scale=0.75]
 				\draw[style=thick,color=blue] (0,0) -- (1,0);
 				\draw[style=thick,color=blue] (1,0) -- (1,1);
 				\draw[style=thick,color=blue] (1,1) -- (0,1);
 				\draw[style=thick,color=blue] (0,1) -- (0,0);
 				\filldraw[fill=black,draw=black] (0,0) circle (0.05);
 				\filldraw[fill=black,draw=black] (0,1) circle (0.05);
 				\filldraw[fill=black,draw=black] (1,0) circle (0.05);
 				\filldraw[fill=black,draw=black] (1,1) circle (0.05);
 			\end{tikzpicture} 
 			\quad $F_9=$
 			\begin{tikzpicture}[baseline=10pt,scale=0.75]
 				\draw[style=thick,color=blue] (0,0) -- (1,0);
 				\draw[style=thick,color=blue] (1,0) -- (1,1);
 				\draw[style=thick,color=blue] (1,1) -- (0,1);
 				\draw[style=thick,color=blue] (0,1) -- (0,0);
 				\draw[style=thick,color=blue] (0,1) -- (1,0);
 				\filldraw[fill=black,draw=black] (0,0) circle (0.05);
 				\filldraw[fill=black,draw=black] (0,1) circle (0.05);
 				\filldraw[fill=black,draw=black] (1,0) circle (0.05);
 				\filldraw[fill=black,draw=black] (1,1) circle (0.05);
 			\end{tikzpicture} 
 			\quad $F_{10}=$
 			\begin{tikzpicture}[baseline=10pt,scale=0.75]
 				\draw[style=thick,color=blue] (0,0) -- (1,0);
 				\draw[style=thick,color=blue] (1,0) -- (1,1);
 				\draw[style=thick,color=blue] (1,1) -- (0,1);
 				\draw[style=thick,color=blue] (0,1) -- (0,0);
 				\draw[style=thick,color=blue] (0,0) -- (1,1);
 				\draw[style=thick,color=blue] (0,1) -- (1,0);
 				\filldraw[fill=black,draw=black] (0,0) circle (0.05);
 				\filldraw[fill=black,draw=black] (0,1) circle (0.05);
 				\filldraw[fill=black,draw=black] (1,0) circle (0.05);
 				\filldraw[fill=black,draw=black] (1,1) circle (0.05);
 			\end{tikzpicture} 
 		
 		\centerline{ Fig. 2: Drawings of $F_i$ $(0\le i\le 10)$ in $\mathcal{F}_4$.}
 		\end{spacing}

 		Let $(G_n)_{n\geq1}$ be a convergent sequence of $F$-free graphs. For simplicity, we write $P(H)=\lim_{n \rightarrow \infty}P(H,G_n)$ and $d(H)=\lim_{n \rightarrow \infty}d(H,G_n)$.  By the definition of $\mathcal{F}_4$, we have the following equality.
 		\begin{equation}\label{EQ: SUMP(F_i)}
 			\sum_{i=0}^{10}P(F_i)=1.
 		\end{equation}
 		By the law of total probability, the (noninduced) density of $C_4$ can be expressed as the sum of the induced densities of graphs on four vertices in the following way:
 		$$d(C_4)=\sum_{i=0}^{10}P(F_i)\cdot N(C_4, F_i).$$
 		Since $N(C_4, F_i)=0$ for $0\le i\le 7$ and $N(C_4, F_8)=N(C_4, F_9)=1$ and $N(C_4, F_{10})=3$, we have 
 		\begin{equation}\label{EQ: d(C_4)}
 			d(C_4)=P(F_8)+P(F_9)+3P(F_{10}).
 		\end{equation}
 		By Lemma~\ref{LEM: K_4}, we have 
 		$$P(F_{10})=P(K_4)=\lim_{n \rightarrow \infty}\frac{N(K_4, G_n)}{{n\choose 4}}\leq \frac{r^3-6r^2+11r-6}{r^3}.$$
So
 		\begin{equation}
 			\begin{aligned}
 				Q_0(r):&=\sum_{i=0}^{9}(\frac{r^3-6r^2+11r-6}{r^3})\cdot P(F_i)+P(F_{10})\cdot\frac{-6r^2+11r-6}{r^3}\\
 				&=\frac{r^3-6r^2+11r-6}{r^3}-P(F_{10})\geq 0\\
 			\end{aligned}
 			\nonumber
 		\end{equation}
 		
In the following computations, we will use two sets of flags $\mathcal{F}^{\sigma_1}_3$ and $\mathcal{F}^{\sigma_2}_3$, where

 		~~~~~~~~~~~~~~~~~~~~~~$\sigma_1=$
 		\begin{tikzpicture}[baseline=0.01pt]
 			\filldraw[fill=yellow,draw=black] (0,0) rectangle (0.1,0.1) node[below left]{\tiny $1$};
 			\filldraw[fill=yellow,draw=black] (1,0) rectangle (1.1,0.1) node[below right]{\tiny $2$};
 		\end{tikzpicture},\quad
 		$\sigma_2=$
 		\begin{tikzpicture}[baseline=0.01pt]
 			\draw[style=thick,color=blue] (0.1,0.05) -- (1.0,0.05);
 			\filldraw[fill=yellow,draw=black] (0,0) rectangle (0.1,0.1) node[below left]{\tiny $1$};
 			\filldraw[fill=yellow,draw=black] (1,0) rectangle (1.1,0.1) node[below right]{\tiny $2$};
 		\end{tikzpicture}.\\
By (\ref{EQ: sigma2}), we have $P\left(\llbracket \alpha^{2} \rrbracket_{\sigma_i}\right) \geq 0$, where  $\alpha \in \mathcal{A}^{\sigma_i}$ and $\alpha^2=\alpha\cdot\alpha$ for $i=1,2$. Therefore,  each of the following three expressions is nonnegative for all $r\geq 3$.\\
$
\begin{aligned}
 			 Q_1(r)&=6P\Big( \big\llbracket \big( (r-1)
 		\begin{tikzpicture}[baseline=8pt,scale=0.75]
 			\filldraw[fill=yellow,draw=black] (0,0) rectangle (0.1,0.1) node[below left]{\tiny 1};
 			\filldraw[fill=yellow,draw=black] (1,0) rectangle (1.1,0.1) node[below right]{\tiny 2};
 			\filldraw[fill=black,draw=black] (0.5,1) circle (0.05);
 		\end{tikzpicture}-
 		\begin{tikzpicture}[baseline=8pt,scale=0.75]
 			\filldraw[fill=yellow,draw=black] (0,0) rectangle (0.1,0.1) node[below left]{\tiny 1};
 			\filldraw[fill=yellow,draw=black] (1,0) rectangle (1.1,0.1) node[below right]{\tiny 2};
 			\draw[style=thick,color=blue] (1.05,0.1) -- (0.5,1.0);
 			\draw[style=thick,color=blue] (0.05,0.1) -- (0.5,1.0);
 			\filldraw[fill=black,draw=black] (0.5,1) circle (0.05);
 		\end{tikzpicture}\big)^2\big\rrbracket_{\sigma_1}\Big)\\
 	&= (6r^2-12r+6)P(F_0)+(r^2-2r+1)P(F_1)+(1-r)P(F_2)\\
 	&\quad +(3-3r)P(F_3)+2P(F_8)+P(F_9)
\end{aligned} 
$\\
$ 	\begin{aligned}
 	 Q_2(r)&=6P\Big( \big\llbracket \big(
 		\begin{tikzpicture}[baseline=8pt,scale=0.75]
 			\filldraw[fill=yellow,draw=black] (0,0) rectangle (0.1,0.1) node[below left]{\tiny 1};
 			\draw[style=thick,color=blue] (0.1,0.05) -- (1.0,0.05);
 			\filldraw[fill=yellow,draw=black] (1,0) rectangle (1.1,0.1) node[below right]{\tiny 2};
 			\draw[style=thick,color=blue] (0.05,0.1) -- (0.5,1.0);
 			\filldraw[fill=black,draw=black] (0.5,1) circle (0.05);
 		\end{tikzpicture}-
 		\begin{tikzpicture}[baseline=8pt,scale=0.75]
 			\filldraw[fill=yellow,draw=black] (0,0) rectangle (0.1,0.1) node[below left]{\tiny 1};
 			\filldraw[fill=yellow,draw=black] (1,0) rectangle (1.1,0.1) node[below right]{\tiny 2};
 			\draw[style=thick,color=blue] (0.1,0.05) -- (1.0,0.05);
 			\draw[style=thick,color=blue] (1.05,0.1) -- (0.5,1.0);
 			\filldraw[fill=black,draw=black] (0.5,1) circle (0.05);
 		\end{tikzpicture}\big)^2\big\rrbracket_{\sigma_2}\Big) \\		
 	&=3P(F_3)+P(F_7)-P(F_6)-4P(F_8) 		
 \end{aligned}
$\\
$
\begin{aligned}		
Q_3(r)&=6P\Big( \big\llbracket \big( (r-2)
 		\begin{tikzpicture}[baseline=8pt,scale=0.75]
 			\filldraw[fill=yellow,draw=black] (0,0) rectangle (0.1,0.1) node[below left]{\tiny 1};
 			\draw[style=thick,color=blue] (0.1,0.05) -- (1.0,0.05);
 			\filldraw[fill=yellow,draw=black] (1,0) rectangle (1.1,0.1) node[below right]{\tiny 2};
 			\draw[style=thick,color=blue] (0.05,0.1) -- (0.5,1.0);
 			\filldraw[fill=black,draw=black] (0.5,1) circle (0.05);
 		\end{tikzpicture}+(r-2)
 		\begin{tikzpicture}[baseline=8pt,scale=0.75]
 			\filldraw[fill=yellow,draw=black] (0,0) rectangle (0.1,0.1) node[below left]{\tiny 1};
 			\filldraw[fill=yellow,draw=black] (1,0) rectangle (1.1,0.1) node[below right]{\tiny 2};
 			\draw[style=thick,color=blue] (0.1,0.05) -- (1.0,0.05);
 			\draw[style=thick,color=blue] (1.05,0.1) -- (0.5,1.0);
 			\filldraw[fill=black,draw=black] (0.5,1) circle (0.05);
 		\end{tikzpicture}-2
 		\begin{tikzpicture}[baseline=8pt,scale=0.75]
 			\filldraw[fill=yellow,draw=black] (0,0) rectangle (0.1,0.1) node[below left]{\tiny 1};
 			\draw[style=thick,color=blue] (0.1,0.05) -- (1.0,0.05);
 			\filldraw[fill=yellow,draw=black] (1,0) rectangle (1.1,0.1) node[below right]{\tiny 2};
 			\draw[style=thick,color=blue] (1.05,0.1) -- (0.5,1.0);
 			\draw[style=thick,color=blue] (0.05,0.1) -- (0.5,1.0);
 			\filldraw[fill=black,draw=black] (0.5,1) circle (0.05);
 		\end{tikzpicture}
 	\big)^2\big\rrbracket_{\sigma_2}\Big)\\
 	&=(3r^2-12r+12)P(F_3)+(r^2-6r+12)P(F_6)+(r^2-8r+12)P(F_7)\\
 	&\quad+(4r^2-16r+16)P(F_8)+(20-8r)P(F_9)+24P(F_{10}).
  \end{aligned}
$\\		
 		In addition, it can be easily checked that for all $r\geq 3$, the following fractions are all nonnegative.\\
$
\begin{aligned}
 q_0(r)&=\frac{3(2r-3)^2}{2(3r^2-11r+9)},\qquad q_1(r)=\frac{2r^2-6r+3}{4r^3(3r^2-11r+9)},\\
q_2(r)&=\frac{8r^2-28r+21}{16(3r^2-11r+9)},\qquad q_3(r)=\frac{8r^2-12r+3}{16r^2(3r^2-11r+9)}.
\end{aligned}	
$\\ 		
So $\sum_{j=0}^{3}q_j(r)Q_j(r)\ge 0$. Therefore, by Equation~(\ref{EQ: d(C_4)}), we have 
 		\begin{eqnarray}
 			d(C_4)&\leq& P(F_8)+P(F_9)+3P(F_{10})+ \sum_{j=0}^{3}q_j(r)Q_j(r)\label{EQ: e5}\\
 			&=& \sum_{i=0}^{10}c_{F_i}P(F_i),\nonumber
 		\end{eqnarray}
 		where  $c_{F_i}$ is the coefficient of $P(F_i)$ after combining like-terms in (\ref{EQ: e5}).
 		The exact values of $c_{F_i}$ are listed in the following.
 		
 		$\bullet$ $ c_{F_0}=c_{F_3}=c_{F_8}=c_{F_9}=c_{F_{10}}=\frac{3(r-1)(r^2-3r+3)}{r^3}.$\\
 		
 		$\bullet$ $c_{F_1}=\frac{26r^5 - 226r^4 + 767r^3 - 1272r^2 + 1029r - 324}{4r^3(3r^2 - 11r + 9)}$\\
 		
 		$\bullet$ $c_{F_2}=\frac{24r^5 - 218r^4 + 758r^3 - 1269r^2 + 1029r - 324}{4r^3(3r^2 - 11r + 9)}$\\
 		
 		$\bullet$ $c_{F_4}=c_{F_5}=\frac{3(2r - 3)^2(r^3 - 6r^2 + 11r - 6)}{r^3(6r^2 - 22r + 18)}$\\
 		
 		$\bullet$ $c_{F_6}=\frac{48r^5 - 448r^4 + 1575r^3 - 2601r^2 + 2070r - 648}{8r^3(3r^2 - 11r + 9)}$\\
 		
 		$\bullet$ $c_{F_7}=\frac{28r^5 - 242r^4 + 804r^3 - 1302r^2 + 1035r - 324}{4r^3(3r^2 - 11r + 9)}$\\
 		By (\ref{EQ: SUMP(F_i)}) and (\ref{EQ: e5}), we have
 		\begin{eqnarray}\label{EQ：upperbound}
 			d(C_4)&\leq& \max\{c_{F_i} : 0\le i\le 10\}\\
 			&=& \frac{3(r-1)(r^2-3r+3)}{r^3},	\nonumber
 		\end{eqnarray}    
 		for all $r\geq 3$, where the equality holds by examining leading coefficients and factoring.
 		This  completes the proof of the claim.                   
 	\end{proof}

 	Let $(G_n)_{n\geq1}$ be a sequence of $F$-free graphs with $\lim_{n \rightarrow \infty} d(C_4, G_n)=OPT_r(C_4)$. Then 
 	\begin{equation}\label{EQ: OPT(C_4)}
 		\lim_{n \rightarrow \infty} d(C_4, G_n)=OPT_r(C_4)=\frac{3(r-1)(r^2-3r+3)}{r^3}\le\lim_{n\rightarrow\infty}\sum_{i=0}^{10}c_{F_i}P(F_i, G_n).
 	\end{equation} 
 	Let $\mathcal{T}=\{F_0, F_3, F_8, F_9, F_{10}\}$, namely,
 	$$\mathcal{T}=\Biggl\{
 	\begin{tikzpicture}[baseline=10pt,scale=0.75]
 		\filldraw[fill=black,draw=black] (0,0) circle (0.05);
 		\filldraw[fill=black,draw=black] (0,1) circle (0.05);
 		\filldraw[fill=black,draw=black] (1,0) circle (0.05);
 		\filldraw[fill=black,draw=black] (1,1) circle (0.05);
 	\end{tikzpicture},\qquad
 	\begin{tikzpicture}[baseline=10pt,scale=0.75]
 		\draw[style=thick,color=blue] (0,0) -- (1,0);
 		\draw[style=thick,color=blue] (0,1) -- (0,0);
 		\draw[style=thick,color=blue] (0,0) -- (1,1);
 		\filldraw[fill=black,draw=black] (0,0) circle (0.05);
 		\filldraw[fill=black,draw=black] (0,1) circle (0.05);
 		\filldraw[fill=black,draw=black] (1,0) circle (0.05);
 		\filldraw[fill=black,draw=black] (1,1) circle (0.05);
 	\end{tikzpicture},\qquad
 	\begin{tikzpicture}[baseline=10pt,scale=0.75]
 		\draw[style=thick,color=blue] (0,0) -- (1,0);
 		\draw[style=thick,color=blue] (1,0) -- (1,1);
 		\draw[style=thick,color=blue] (1,1) -- (0,1);
 		\draw[style=thick,color=blue] (0,1) -- (0,0);
 		\filldraw[fill=black,draw=black] (0,0) circle (0.05);
 		\filldraw[fill=black,draw=black] (0,1) circle (0.05);
 		\filldraw[fill=black,draw=black] (1,0) circle (0.05);
 		\filldraw[fill=black,draw=black] (1,1) circle (0.05);
 	\end{tikzpicture} ,\qquad
 	\begin{tikzpicture}[baseline=10pt,scale=0.75]
 		\draw[style=thick,color=blue] (0,0) -- (1,0);
 		\draw[style=thick,color=blue] (1,0) -- (1,1);
 		\draw[style=thick,color=blue] (1,1) -- (0,1);
 		\draw[style=thick,color=blue] (0,1) -- (0,0);
 		\draw[style=thick,color=blue] (0,1) -- (1,0);
 		\filldraw[fill=black,draw=black] (0,0) circle (0.05);
 		\filldraw[fill=black,draw=black] (0,1) circle (0.05);
 		\filldraw[fill=black,draw=black] (1,0) circle (0.05);
 		\filldraw[fill=black,draw=black] (1,1) circle (0.05);
 	\end{tikzpicture},\qquad
 	\begin{tikzpicture}[baseline=10pt,scale=0.75]
 		\draw[style=thick,color=blue] (0,0) -- (1,0);
 		\draw[style=thick,color=blue] (1,0) -- (1,1);
 		\draw[style=thick,color=blue] (1,1) -- (0,1);
 		\draw[style=thick,color=blue] (0,1) -- (0,0);
 		\draw[style=thick,color=blue] (0,0) -- (1,1);
 		\draw[style=thick,color=blue] (0,1) -- (1,0);
 		\filldraw[fill=black,draw=black] (0,0) circle (0.05);
 		\filldraw[fill=black,draw=black] (0,1) circle (0.05);
 		\filldraw[fill=black,draw=black] (1,0) circle (0.05);
 		\filldraw[fill=black,draw=black] (1,1) circle (0.05);
 	\end{tikzpicture} 
 	\Biggr\}.$$
 	Then $c_{F_i}=\frac{3(r-1)(r^2-3r+3)}{r^3}$ for each $F_i\in\mathcal{T}$. 
 	By (\ref{EQ：upperbound}) and (\ref{EQ: OPT(C_4)}),    for every $F\in \mathcal{F}_4$ with $P(F)>0$, 
 	we have $c_{F}=\frac{3(r-1)(r^2-3r+3)}{r^3}$, which implies that $F\in\mathcal{T}$,  and $P(F)=0$ for all $F\in \mathcal{F}_4\setminus\mathcal{T}$.
 	Notice that none of the graphs in $T$ contain the cocherry graph $P_3^c$ as an induced subgraph. Therefore, we have
 	$$\lim_{n \rightarrow \infty}d(P_3^c, G_n)=\sum_{i=0}^{10}N(P_3^c, F_i)p(F_i)=0.$$
 	This immediately implies the result.     
 \end{proof}    


\noindent{\bf Acknowledgment:}
We thank Professor D\'aniel Gerbner for remarks on the history of the problem of $F$-Tu\'an-good and valuable discussions on this manuscript. 

\noindent\textbf{Data Availability:}
Data sharing is not applicable to this article as no datasets were generated or analyzed during the current study.

\end{document}